\documentclass[11pt,letterpaper,reqno]{amsart}

\makeatletter%
\@mparswitchfalse%
\makeatother%

\reversemarginpar%

\usepackage{amsfonts}
\usepackage{amssymb}
\usepackage{amsmath}
\usepackage{amsthm}
\usepackage{mathtools}
\usepackage{amstext}
\usepackage{tikz-cd}

\usepackage[left=1.3in,right=1.3in,top=1.2in, bottom=1.3in]{geometry}

\definecolor{refkey}{HTML}{0003d6}
\definecolor{labelkey}{HTML}{083b19}

\usepackage[colorinlistoftodos, textwidth=28mm, color=green!30]{todonotes}

\usepackage[ruled,noend,linesnumbered,algosection]{algorithm2e}
\usepackage{graphicx}

\definecolor{cit}{HTML}{117733}
\definecolor{lin}{HTML}{0003d6}
\definecolor{gree}{HTML}{6dc066}
\usepackage{hyperref}
\hypersetup{bookmarksnumbered,colorlinks=true,linkcolor=lin,citecolor=cit,urlcolor=black}

\usepackage[section]{placeins}
\usepackage{caption}

\newcommand{\mycomment}[1]{} %
\newcommand{\mc}[1]{} %

\DeclareMathOperator{\Pic}{Pic}
\DeclareMathOperator{\Cl}{Cl}

\DeclareMathOperator{\CDiv}{CDiv}

\DeclareMathOperator{\rk}{rank}

\newcommand{\Z}{\mathbb{Z}}

\newcommand{\R}{\mathbb{R}}
\newcommand{\PP}{\mathbb{P}}

\newcommand{\Q}{\mathbb{Q}}

\newcommand{\gMC}{generalized Mukai Conjecture}

\DeclareMathOperator{\adjt}{adj}

\DeclareMathOperator{\AD}{AD}
\DeclareMathOperator{\conv}{conv}

\DeclareMathOperator{\relint}{relint}
\DeclareMathOperator{\Int}{int}

\newcommand{\vv}[1]{\lvert#1\rvert}
\newcommand{\vr}[1]{\langle#1\rangle}
\newcommand{\ra}{\to}

\newcommand{\txt}[1]{\text{ #1}} %
\newcommand{\gptversion}{6.1 Sol}

\theoremstyle{plain}%
\newtheorem{theorem}{Theorem}[section]

\newtheorem{corollary}[theorem]{Corollary}
\newtheorem{prop}[theorem]{Proposition}
\newtheorem{lemma}[theorem]{Lemma}

\theoremstyle{definition}%
\newtheorem{defi}[theorem]{Definition}

\newtheorem{example}[theorem]{Example}

\newcommand{\pf}{{\em Proof}}%
\newcommand{\pfof}[1]{{\em Proof of #1}}%

\begin{document}
\title[Gorenstein toric Fano $d$-varieties have Picard ranks at most $2d$]{Gorenstein toric Fano $d$-varieties have Picard ranks at most $2d$}
\author{Zhuang He \and Benjamin Nill}
\email{he.zhu@northeastern.edu}
\email{benjamin.nill@ovgu.de}

\begin{abstract}
We show that the Picard rank of $d$-dimensional Gorenstein toric Fano varieties is at most $2d$. This generalizes Batyrev's conjecture of Picard ranks for smooth toric Fano varieties, which was proved by Casagrande for all $\Q$-factorial Gorenstein toric Fano varieties. The convex-geometric proof is a generalization of Casagrande's proof and Eikelberg's theory of non-$\Q$-factorial toric Picard groups, with key ideas from work of Debarre, {\O}bro, Fujita and the first author. The proof was found using ChatGPT 6.1.
\end{abstract}
\maketitle

\section{Introduction}\label{Intro}
Smooth Fano varieties are among the most studied objects in algebraic geometry. The classification and study of topological invariants of Fano varieties have significance in birational geometry and mirror symmetry. Smooth Fano varieties of dimension $d$ form bounded families under deformation \cite{KMM}. Since the Picard ranks $\rho$ and the second Betti number coincide for smooth Fano varieties, they are invariant within each deformation family, and there can only be finitely many Picard ranks in each dimension. In low dimensions, the maximal Picard ranks are $9$ for smooth Fano surfaces, $10$ for smooth Fano threefolds \cite{MoriMukai1981}, and $18$ for smooth Fano fourfolds \cite{Casagrande2024}. The question of bounding the Picard ranks is open in higher dimensions.

Toric Fano varieties, in particular, Gorenstein toric Fano varieties, are essential in Batyrev's theory \cite{Batyrev1994} of combinatorial mirror symmetry. By the toric dictionary, isomorphism classes of $d$-dimensional Gorenstein toric Fano varieties correspond one-to-one to isomorphism classes\footnote{Here we use the notion of {\em isomorphism} which is also called {\em unimodular equivalence}, namely, given by $x \mapsto Ax +b$ for $A \in \text{GL}(d,\Z)$ and $b \in \Z^d$.} of $d$-dimensional {\em reflexive} polytopes: lattice polytopes whose polar is also a lattice polytope. Arbitrary toric Fano varieties can have Picard ranks as large as possible in all dimensions starting from two. However, there are only finitely many reflexive polytopes in each dimension up to isomorphism, giving finitely many possible Picard ranks. Batyrev has conjectured that the Picard ranks for smooth Fano toric varieties are at most $2d$. Multiple estimates had been found in the smooth case \cite{Voskresenskiĭ1985,Debarre2000}, and the $\Q$-factorial case \cite{Nill2005M}, which is more general, before Casagrande \cite{Casagrande2005} proved the conjecture for all $\Q$-factorial Gorenstein toric Fano varieties. Further simplifications, classifications and computations include \cite{2008Øbro,NillObro2010LargePicardNumber,AJP2014,AssarfNill2016Splitting}. In particular, the equality holds exactly for the product $(S_3)^{d/2}$, where $S_3$ is the del Pezzo surface of degree $6$, isomorphic to the blow-up of $\PP^2$ along three non-collinear points. The corresponding reflexive polytope is the direct sum of $d/2$ copies of the reflexive hexagon $\mathcal{H}$ with vertices $\{(\pm 1, 0), (0, \pm 1), \pm (1,1)\}$.

In the non-$\Q$-factorial case,  or equivalently, for non-simplicial reflexive polytopes, the question has remained open. Let $X$ be a projective Fano variety. The {\em pseudo-index} of $X$ is defined as $\iota_X:=\min \{-K_X \cdot C : C \txt{is rational curve on } X\}$. Introduced in \cite{BCDD2003}, the {\em\gMC{}} claims that if $X$ is a smooth Fano variety, then
\[
\rho_X (\iota_X -1 )\le d.
\]
Fujita \cite{Fujita2019} proved the \gMC{} for toric log Fano pairs, including the Gorenstein ones. Despite this, in any dimension $d\ge 3$ there are non-$\Q$-factorial Gorenstein toric Fano varieties with pseudo-index one (for example, $X\times (\PP^1)^{d-3}$, where $X$ is the toric threefold of the cube $[-1,1]^3$), in which case the \gMC{} holds emptily and provides no restriction on the Picard rank. In dimension three, an exhaustive computation over the $4,319$ three-dimensional reflexive polytopes returns Picard ranks of at most $6$, using Kreuzer and Skarke's database \cite{Kreuzer1998}.

The Picard rank $\rho_X$ of a toric variety $X=X_P$ from the spanning fan over a lattice polytope $P$ is also bounded by $\rk \Cl(X_P)=\vv{V(P)}-d$, which is strictly larger than $\rho_X$ when $X$ is non-$\Q$-factorial. Some bounds of the vertex count and the Picard rank that depend on the specific polytopes were given by \cite{Nill2005M, He2026}. For an example of a simple reflexive polytope with an exponential number of vertices $\vv{V(P)} =6^{d/2}$, just consider $P$ as the direct product $\mathcal{H}^{d/2}$. In the recent preprint \cite{HofscheierNill2026ReflexivePolytopes} it was surprisingly shown that the number of vertices of a reflexive polytope can even grow doubly-exponentially in the dimension $d$. In stark contrast, there is a not enough known result by Eikelberg \cite[Corollary 2.18]{Eikelberg1993} that simple polytopes of dimension at least three give toric varieties of Picard rank one - no matter how large their number of vertices is! This inspired several ideas in the later paper \cite{He2026} by the first author and is one of the chief reasons why we believed that a linear upper bound on the Picard number should exist for all reflexive polytopes. In the present paper, we provide such a sharp upper bound on the Picard number.

\begin{theorem}\label{GPT}
Let $X$ be a Gorenstein toric Fano variety of dimension $d$, then $\rho_X\le 2d$. The equality holds if and only if $X\cong (S_3)^{d/2}$. In particular, if $d$ is odd, then $\rho_X\le 2d-1$. 
\end{theorem}

We found the proof of Theorem \ref{GPT} by using ChatGPT \gptversion{}. Although the proof is AI-generated, the upper bound $2d$ itself was not: we set the goal of finding a linear upper bound, and gave the sharp bound $2d$ as part of the prompts to ChatGPT. The idea that the non-$\Q$-factorial Picard ranks can be no bigger than $2d$ was introduced to the first author by Casagrande in 2024, noticing the computational evidence in dimension $3$, and Fujita's proof \cite{Fujita2019} of the \gMC{} in the non-$\Q$-factorial case. Failed attempts in dimension four and higher to construct reflexive polytopes of large Picard ranks, and the upper bounds found in \cite{He2026} further convinced us that a linear upper bound may be within reach.

The proof of Theorem \ref{GPT} is a synthesis of multiple ideas in literature. Let us give a short overview. Casagrande's proof \cite{Casagrande2005} and the later simplification by \O bro \cite{2008Øbro} that simplicial reflexive polytopes of dimension $d$ have at most $3d$ vertices are based on the following observation: fixing any facet $F$ of a reflexive polytope $P$, the unique primitive inner normal vector $u_F$ of $F$ is a lattice vector, and hence $\vr{u_F,v}\in \Z_{\ge -1}$ for any vertex of $P$. Thus we can group the vertices of $P$ by their `levels' with respect to $u_F$: those on level $-1$ are exactly the vertices of $F$; if $v$ has level zero, then one shows that $v$ must be contained in a facet adjacent to $F$, hence up to $d$ many; otherwise, if $v$ has positive level, then one can choose a so-called `special facet' $F$ such that $D:=\sum_{V\in V(P)} v$ falls in the cone over $F$. As a result, the weighted level-sum of the positive-leveled vertices is at most the level-sum of level-$(-1)$ vertices, which equals $d$. This means that there are also at most $d$ positive-leveled vertices.

In the non-simplicial case, the estimate above is broken, as when $F$ is not a simplex, both the number of vertices and adjacent facets of $F$ can grow fast, leaving no linear upper bounds for any of the three parts. To estimate the Picard rank instead of the vertex count, we need the description by \cite{Eikelberg1992} of the non-$\Q$-factorial toric Picard group, that the space of torus-invariant $\Q$-Cartier divisors on $X$ is dual to the quotient space $H_\Q/\AD(X)$, where $H_\Q$ is the free $\Q$-vector space on the vertices of $P$, and $\AD(X)$ is the space of affine relations that are linear on each face cone (see Section \ref{Prelim}). A {\em numerical basis}, by \cite{Fujita2019} (Definition \ref{NB}), is a subset of vertices $B\subseteq V(P)$ whose classes form a basis of $H_\Q/\AD(X)$. As a result, $\rho_X=\vv{B}-d$.

The proof of Theorem \ref{GPT} starts with choosing a special facet $F$ of $P$, and a numerical basis $B$ that minimizes the level-sum $F(B):=\sum_{v\in B} \vr{u_F, v}$. Such $F$ and $B$ exist, as both are chosen in a finite list. Then, as in the simplicial case, we subdivide $B$ by levels. Using results on numerical bases and studying more closely how to cover the vertices on level $0$ we derive the following estimates:
\[|B_{-}|=d,\quad |B_0| \le 2d-m, \quad |B_{+}| \le m,\]
where $m$ is given in Proposition~\ref{zeroCovering}. This gives $|B|\le 3d$, so $\rho_X \le 2d$, as desired.

We will discuss and analyze the proof and the origin of the ideas in more detail in Section~\ref{AI-subsec}.

{\bfseries Declaration of AI-generated contents} The paper is written by human authors with AI-generated results. The proof of Theorem \ref{GPT} is generated by ChatGPT \gptversion{}. We verified the AI-generated proof, and then rewrote it with our interpretation, and moderate simplification and consolidation. The authors take full responsibility for the correctness of the proof, and the accuracy of credits and references in the paper. 

{\bfseries Acknowledgment} 
We thank Cinzia Casagrande for communications and suggestions. 
The first author thanks the Otto-von-Guericke-Universität Magdeburg for hosting the research visit, during which most of the paper was developed. 
The first author was partially funded by the European Union - Next Generation EU, Mission 4, Component 1, CUP D53D23005860006, within the project PRIN 2022L34E7W ``Moduli spaces and birational geometry''. The second author is funded by the Deutsche Forschungsgemeinschaft (DFG, German Research Foundation) – 539867500 as part of the research priority program Combinatorial Synergies. 

\section{Preliminaries}\label{Prelim}

\subsection{Polyhedral geometry and reflexive polytopes} Consider $N=\Z^d$. Let $N_\R:=N\otimes_\Z\R$. For any point set $S\subseteq N_\R$, $\conv(S)$ is the {\em convex hull} of $S$. The interior and relative interior of a set $A\subseteq N_\R$ is denoted by $\Int A$ and $\relint A$ respectively. A $d$-dimensional polytope in $N_\R$ is referred to as a {\em $d$-tope}. Let $P$ be a $d$-tope. A {\em facet} of $P$ is a $(d-1)$-face, and a {\em ridge} of $P$ is a $(d-2)$-face. Let $F(P)$ and $V(P)$ be the set of facets and vertices of $P$. 

Let $M=N^*$ be the dual lattice. A full dimensional lattice polytope $P\subset N_\R$ is called {\em reflexive} if $0\in \Int P$, and for every facet $F$ of $P$,
\[ \vr{u_F, x}=-1 \txt{for all } x \in F,\]
where $u_F$ is the primitive inner normal vector of $F$ (i.e., $u_F\in M$ and $\lambda u_F \not\in M$ for all $0 < \lambda <1$). Recall the definition of the {\em polar} or {\em dual} polytope:
\[P^*:=\{y \in M_\R \;:\; \vr{y,x} \ge -1 \txt{for all } x \in P\}.\]
A lattice polytope $P$ containing $0$ in its interior is reflexive if and only if $P^*$ is a lattice polytope too. In this case, the set of vertices of $P^*$ are precisely the primitive inner normal vectors of facets of $P$. 

Two distinct facets $F$ and $G$ of $P$ are {\em adjacent} if they share the same ridge. Write $\adjt(F)$ as the set of adjacent facets of $F$. %
For a facet $F$ of $P$, let $Z_F$ be the parallel hyperplane in $N_\R$ defined by 
\[Z_F := \{x \in N_\R:\vr{u_F, x}=0\}.\] If $P$ is reflexive, then by definition, $\vr{u_F, x}\in \Z_{\ge -1}$, and for any vertex $x$ of $P$ we have either $x\in F$, or $x\in Z_F$, or otherwise, $\vr{u_F, x}\in \Z_{\ge 1}$.	

We denote by $X=X_P$ the Gorenstein toric Fano variety associated to a reflexive polytope $P$, where the fan of $X$ is the spanning fan $\Sigma(P)$ over $P$. We refer to \cite{Nill2005} for more on reflexive polytopes and their associated toric varieties.

\subsection{Picard groups and numerical bases} We follow \cite[\S 3.1]{He2026} and \cite[\S 3]{Fujita2019}.

Eikelberg \cite{Eikelberg1992} described the Picard group of projective toric varieties under even more general assumptions than needed here. Let us recall his description.

Let $P$ be a reflexive $d$-tope, $n:=\vv{V(P)}$, $H:= \Z^{\oplus \vv{V(P)}}\cong \Z^n$, with the standard basis $h_1,\cdots, h_n$ corresponding to the vertices $\{v_i\}_{i=1}^n$ of $P$. Let $X=X_P$. For each facet $F$ of $P$, define:
\begin{equation}
\AD(F):=\left\{\sum_{i=1}^n a_i h_i\in H_\Q: 
\begin{array}{l}
a_i =0 \txt{if } v_i \not\in F,\\ \sum_{i=1}^n a_i v_i= 0, \txt{and } \sum_{i=1}^n a_i =0.
\end{array}
\right\} \subseteq H_\Q.
\end{equation}

Define $\AD(X):= \sum_{F \in F(P)} \AD(F) \subseteq H_\Q$. Let $\CDiv_T(X)$ denote the set of torus-invariant Cartier divisors on $X$, and $D_1, \ldots, D_n$ the prime divisors corresponding to the $n$ rays of the spanning fan of $P$. A crucial result is:

\begin{lemma}\cite[Theorem 3.2]{Eikelberg1992}\label{EikelbergT}
There exists a perfect pairing $\CDiv_T(X)_\Q \times H_\Q/\AD(X)$ $\ra \Q$, sending $(\sum_i r_i D_i, \sum_i a_i h_i)$ to $\sum_i r_i a_i$. The Picard group $\Pic(X)$ is a free abelian group of rank $n-\dim \AD(X)-d$.
\end{lemma}

Hence, the {\em Picard number} or {\em Picard rank} of $X$ equals
\[\rho_X = n-\dim \AD(X)-d.\]

\begin{defi}\cite[Definition 3.5]{Fujita2019}\label{NB}
Given $X=X_P$ as above. We say that a set of vertices $B=(v_1,\cdots, v_{d+\rho_X})$ is a {\em numerical basis} of $X$ (or $P$), if the classes of $h_1, \cdots, h_{d+\rho_X}$ in $H_\Q/\AD(X)$ form a basis of $H_\Q/\AD(X)$.
\end{defi}

\begin{lemma}[See the text following \cite{Fujita2019}, Definition 3.5]\label{NBexists}
There exists a numerical basis $B$.
\end{lemma}
\pf. The quotient map $H_\Q\ra H_\Q/\AD(X)$ is surjective, and $\{h_1,\cdots, h_n\}$ span $H$, so their classes span $H_\Q/\AD(X)$ too. Thus they must contain a basis of $H_\Q/\AD(X)$.\qed

By definition, all numerical bases of $X$ are of the same length $d+\rho_X$. We also observe that $F$ is simplex if and only if $\AD(F)=0$. In particular, if $P$ is simplicial, then $\AD(X)=0$, so $V(P)$ gives the unique numerical basis. 

Let us describe some properties of numerical bases that will be essential in the proof of the main result.

\begin{lemma}\cite[Prop. 3.6]{Fujita2019}\label{negativeCount}
If $B$ is a numerical basis of $X$, then for any facet $F$ of $P$, $\vv{B\cap V(F)}\le d$.
\end{lemma}

For any $A\subseteq V(P)$, let $H(A)$ be the subspace of $H_\Q/\AD(X)$ spanned by the quotient classes $\{\overline{h}_i: v_i \in A\}$.

\begin{lemma}\label{mostD}
Let $G$ be a $k$-dimensional face of $P$, with $k<d$. Then $\dim H(V(G)) = k+1$.
\end{lemma}
\pf. This follows from the slightly more general result \cite[Lemma 3.10]{He2026}.\qed

\begin{lemma}\label{1SF}
Let $B$ be a numerical basis of $X$. Let $F$ be a facet of $P$, and $S\subseteq V(P)$. If $S\cap F$ contains all the vertices of $F$ on a ridge $\tau$ of $F$, then $\dim H(S\cup V(F))\le \dim H(S)+1$.
\end{lemma}

Lemma \ref{1SF} is essentially proved in \cite[Lemma 2.16]{Eikelberg1993}, and is the key step of \cite[Thm 3.19]{He2026}.

\pfof{Lemma \ref{1SF}}. The dimension formula yields
\begin{align*}
\dim H(S\cup V(F)) &= \dim H(S) + \dim H(V(F)) - \dim H(S)\cap H(V(F))\\
									&\le \dim H(S) + \dim H(V(F))- \dim H(S\cap V(F)).
\end{align*}

By Lemma \ref{mostD}, $ \dim H(V(F))=d$, and $\dim H(S\cap V(F))\ge \dim H(V(\tau))= d-1$, so $\dim H(S\cup V(F))\le \dim H(S)+1$. \qed

\begin{corollary}\label{1perstep}
Let $B$ be a numerical basis of $X$. Let $G_0, \cdots, G_k$ be distinct facets of $P$ such that for every $i=1, \ldots, k$, $G_i$ is adjacent to a facet $G_{j_i}$, with $0\le j_i< i$. Then 
\[\vv{B\cap (\cup_{i=0}^k V(G_i))}\le d+k.\]
\end{corollary}

\begin{example}\label{1374}
We consider the reflexive $3$-tope $P$ in \cite[Example 3.6]{He2026}. As shown in Figure \ref{1374figure}, $\vv{V(P)}=7$, and $\vv{F(P)}=8$. $\AD(X)=\Q(-h_2+h_3+h_5-h_6, -h_0+h_1+h_5-h_6)$, so $\rho_X= 7-2-3=2$. A five-vertex subset $B$ is a numerical basis if and only if the two vertices that $B$ misses are from two different subsets out of $\{h_0,h_1\},\{h_2,h_3\}$ and $\{h_5, h_6\}$. Thus $P$ has exactly $12$ numerical bases. Note that we can cover all the vertices of $P$ by the facet sequence $(\sigma_{2356}, \sigma_{0156}, \sigma_{014})$. Thus by Corollary \ref{1perstep}, $ \dim H_\Q/\AD(X) \le 3+2=5$, so we get $\rho_X\le  \dim H_\Q/\AD(X) - 3 = 2$.
\begin{figure}[h]
\includegraphics[width=0.25\textwidth]{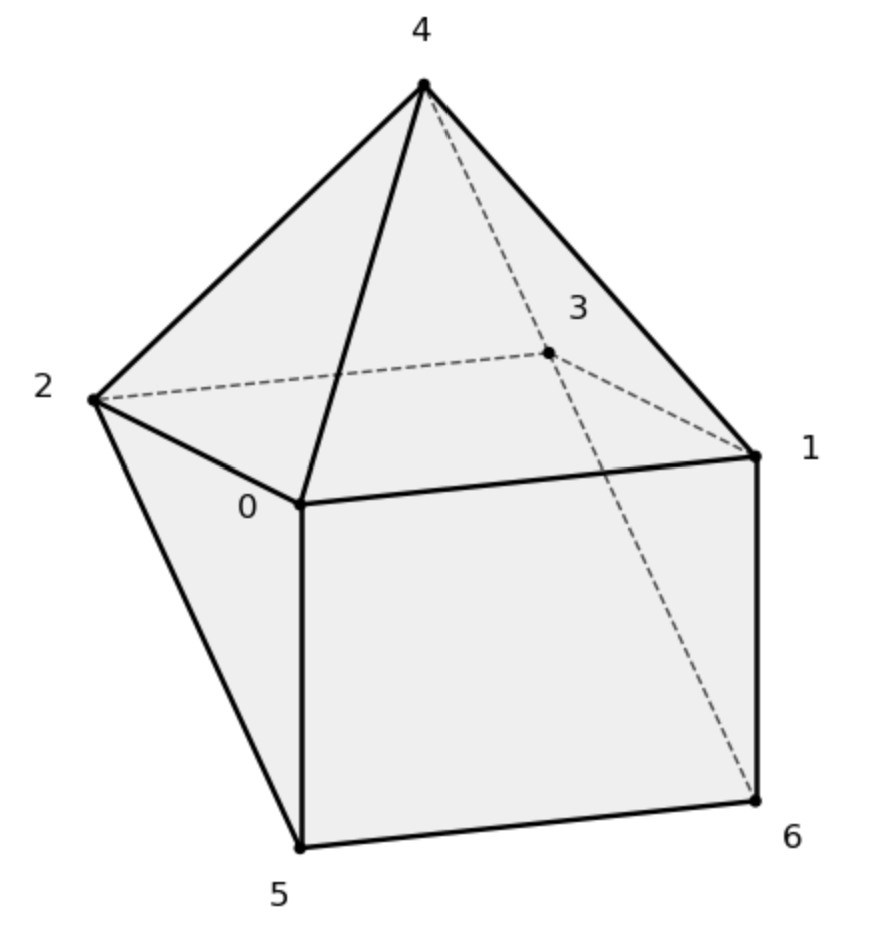}
\caption{A reflexive $3$-tope with Picard rank $2$. Each quadrilateral facet decides a generator of $\AD(X)$.}\label{1374figure}
\end{figure}
\end{example}

\section{Linear upper bound of Picard ranks}\label{GPTproof}
\subsection{Covering the vertices on the zero layer}

A key ingredient of Casagrande's proof \cite{Casagrande2005} of the vertex-count upper bound in the simplicial case is the description of zero-level vertices. In \cite[Lemma 5.5]{Nill2005M} it was shown by the second author that for any facet $F$ of a reflexive polytope $P$ the following holds: any vertex of $P$ in $Z_F$ (so, on the parallel hyperplane through the origin) is contained in a facet of $P$ that is adjacent to $F$. If $P$ is simplicial, then any facet $F$ has exactly $d$ neighboring facets, so any $Z_F$ contains at most $d$ vertices. This is no longer the case when $P$ is non-simplicial. However, as the following result shows we can still cover the vertices of $Z_F$ by at most $2d-2$ many adjacent facets.

\begin{prop}\label{zeroCovering}
Let $P$ be a reflexive $d$-tope with $d\ge 2$. Let $F$ be a facet of $P$. The following holds:

\begin{enumerate}
  \item There are $m$ facets $G_1,\cdots, G_{m}$ in $\adjt(F)$ (with $2\le m\le d$), such that for some $\lambda_1\, \ldots, \lambda_m \in \Q_{>0}$ and $r \in \Q$ we have 
\[\sum_{i=1}^{m} \lambda_i u_{G_i} = r u_F,\]
and $u_{G_1}, \ldots, u_{G_m}$ restricted to $Z_F$ are affinely independent.
\item For every choice of $G_1, \ldots, G_m$ in (1) there exist at most $2d-m$ facets in $\adjt(F)$, containing $G_1, \ldots, G_m$, which cover all the vertices of $P$ on $Z_F$.
\item For any facet $G \in \adjt(F)$ there exist $m$ facets $G_1,\cdots, G_m$ satisfying (1) with $G_1=G$.
\end{enumerate}
\end{prop}

Proposition \ref{zeroCovering} is proved by a refinement of Steinitz’s theorem on convex hulls. There are multiple ways stating the theorem, and we use the following form.

\begin{lemma}\label{SteinitzRefined}
Let $S$ be a set of points in $\R^l$. Suppose $x\in \Int\conv(S)$, and there exist affinely independent points $T$ in $S$ such that $x\in \relint\conv(T)$ and $1\le \dim(T)\le l$. Then there exists $S'\subseteq S$ such that $T\subseteq S'$, $\vv{S'}\le 2l-\dim(T)+1$, and $x\in\Int\conv(S')$.
\end{lemma}

\pfof{Lemma \ref{SteinitzRefined}}. This follows directly from \cite[Lemma 3.4]{CJT2026}, by assigning their $(X,Z,Y\cup Z,k,k')$ to our $(S,T,S',l,\dim(T))$ respectively.\qed

 \pfof{Proposition \ref{zeroCovering}}.  
Let us use for convenience $N=\Z^d$ and $M=(\Z^d)^*$.

Ad (1): We define the projection $\pi$ of $(\R^d)^*$ along the vertex $u_F$ of $P^*$. Note that its image is the dual space to $Z_F$. By considering the tangent cone of the vertex $u_F$ of $P^*$, it is clear that the convex hull of $S := \{\pi(u_G)\}_{G \in \adjt(F)}$ contains the origin in its interior. Choose a minimal subset $T \subseteq S$ that contains the origin in the relative interior of its convex hull, clearly, $T$ is affinely independent. Define $m := |T|$. We observe $2 \le m \le \dim(T) + 1 \le (d-1)+1=d$. Let $G_1, \ldots, G_m \in \adjt(F)$ such that $T = \{\pi(u_{G_i}) \;:\; i=1, \ldots, m\}$. From this the statement follows.

Ad (2): We apply Lemma~\ref{SteinitzRefined} to $S$ and $T$ (with $l=d-1$ and $\dim(T)=m-1$) and get the existence of $S' \subseteq S$ such that $T \subseteq S'$, $|S'| \le 2(d-1)-(m-1)+1=2d-m$, and $0 \in \Int\conv(S')$. Let $v$ a vertex of $P$ contained in $Z_F$. As $\conv(S')$ spans the dual space of $Z_F$ and $0 \in \Int\conv(S')$, there must exist some element $\pi(u_H)$ in $S'$ (for some $H \in \adjt(F)$) with $\vr{u_H,v}=\vr{\pi(u_H),v} < 0$. As $P$ is reflexive, so $\vr{u_H,v} \in \Z$, and $u_H \in P^*$, this implies $\vr{u_H,v} = -1$, so $v \in H$, as desired. 

Ad (3): Since $0\in \Int \conv(S)$, the cones generated by $S$ cover the whole space $\pi(M_\R)$. Therefore $-\pi(u_G)$ is a nonnegative sum of some $U\subseteq S\backslash \{\pi(u_G)\}$. Choose a minimal generating set $T'$ from $U$, then $T'$ is linearly independent, and generate $-\pi(u_G)$ with positive coefficients. Let $T:=\{\pi(u_G)\} \cup T'$. Then $\conv(T)$ contains the origin in the relative interior, and $T$ is affinely independent, thus (1) follows.\qed

\subsection{Proof of the main result}\

\pfof{Theorem \ref{GPT}}. The claim is true for $d=1$, as then $X\cong \PP^1$ with Picard rank one. Assume then $d\ge 2$. Let $\mathcal{B}$ be the set of numerical basis on $P$. Then $\mathcal{B}$ is a finite set. Therefore there exists a numerical basis $B\in \mathcal{B}$ and a facet $F$ of $P$ minimizing the following integer:
\begin{align}
F(B):=\sum_{b\in B} \vr{u_F, b}
\end{align}
among all the $B'\in \mathcal{B}$ and $F'\in F(P)$. We say a pair $(B',F')$ is minimizing if $F'(B')=F(B)$ achieves the minimum. Let $B_{-}$, $B_0$ and $B_{+}$ be the vertices in $B$ such that $\vr{u_F, w}$ is negative, zero and positive respectively. Note $B_{-}\subset V(F)$. %

By Proposition \ref{zeroCovering} (1), there exists $m$ distinct facets $G_1,\cdots, G_{m}$ neighboring $F$, such that $\sum_{i=1}^{m} \lambda_i u_{G_i} = r u_F$, for $\lambda_i\in\Q_{>0}$ and $r\in\Q$. Then by Proposition \ref{zeroCovering} (2), there are at most $2d-m$ distinct facets of $P$ adjacent to $F$ that cover all the vertices of $P$ on $Z_F$, and hence they also cover $B_0$. By Corollary~\ref{1perstep}, $\vv{B_{-}}+\vv{B_0}\le d+2d-m=3d-m$.

We now show that $\vv{B_+}\le m$, so together we have 
\begin{align}\label{3parts}
\vv{B}= \vv{B_{-}}+\vv{B_0}+\vv{B_{+}}\le 3d-m+m=3d,
\end{align}
and then $\rho_X=\vv{B}-d\le 2d$, which proves the inequality of the theorem. 

{\bfseries Proof of $\vv{B_+}\le m$}:
Firstly, we show $\vv{B_{-}}=d$. Suppose  $\vv{B_{-}}<d$. By Lemma \ref{mostD}, $\dim H(V(F)) = d$. Thus there is a vertex $w\in V(F)$ such that the generator $h_w$ it indexes is not spanned by $\{\overline{h}_v: v\in B_{-}\}$ in $H_\Q/\AD(X)$. Write $h_w$ as a sum under the basis $\{h_b\}_{b\in B}$. Then there is a vertex $y\in B\backslash B_{-}$ such that the coefficient of $h_w$ at $h_b$ is nonzero. Let $B':=(B\backslash \{y\})\cup \{w\}$. Then $\{h_v: v\in B'\}$ also generate $H_\Q/\AD(X)$, so $B'$ is a numerical basis too. Now $F(B')<F(B)$ because $F(B')-F(B)=\vr{u_F, w-y}= -1-\vr{u_F, y}$, and $\vr{u_F, y}\ge 0$. This contradicts that $B$ is minimizing.

For each $i$, $1\le i\le m$, we construct another numerical basis $B_i$ as follows. Since $G_i$ is adjacent to $F$, $G_i$ share a ridge $\sigma_i$ with $F$. Choose $d-1$ affinely independent vertices $E_i$ of $\sigma_i$, and choose any other vertex $z_i$ in $V(F)\backslash V(\sigma_i)$. Then let $B_i:=E_i\sqcup \{z_i\} \sqcup B_0 \sqcup B_{+}$. Since $E_i\sqcup \{z_i\}$ also span the subspace $H(V(F))$, we find $B_i$ span $H_\Q/\AD(X)$ too, so the length count shows that $B_i$ is also a numerical basis of $P$. 

Write $u_i$ for $u_{G_i}$. Let $w:=w(B, F):= \sum_{b\in B_0\cup B_{+}} b$, and $a :=\vr{u_F, w}$. Then $F(B)=-d + a$. Note that $|B_+| \le \vr{u_F,\sum_{b \in B_{+}} b} = a$. In particular, it remains to show $a \le m$. We have
\[
G_i(B_i)=\sum_{b\in B_i} \vr{u_{i}, b}=-(d-1) + \vr{u_i, z_i} + \vr{u_i, w}.
\]
The term $d-1$ is from the fact that $P$ is reflexive, so $\vr{u_i, x}=-1$ for any $x\in G_i$. By the choice, $(B,F)$ is minimizing, so $F(B)\le G_i(B_i)$. Therefore
\begin{align}\label{x1}
a\le 1 + \vr{u_i, z_i} + \vr{u_i, w}.
\end{align}
Sum up (\ref{x1}) with weights $\lambda_i$, for $i$ from $1$ to $m$:
\[ a\sum_{i=1}^m \lambda_i \le \sum_{i=1}^m\lambda_i(1+\vr{ u_i, z_i})+ \left\langle\sum_{i=1}^m\lambda_i u_i, w\right\rangle.\]
\[ a\left(\sum_{i=1}^m \lambda_i-r\right) \le \sum_{i=1}^m\lambda_i(1+\vr{ u_i, z_i}).\]
On the other hand, for each fixed $k$, $\sum_j \lambda_{j=1}^m \vr{u_j, z_k}= r\vr{u_F, z_k}=-r$, as $z_k\in F$. So
\begin{align}\label{x2}
\sum_{i=1}^m \lambda_i-r= \sum_{j=1}^m \lambda_j (1+ \vr{u_j, z_k}). 
\end{align}
Since $z_k\in P$, $\vr{u_j, z_k}\ge -1$. Since $z_k\not\in G_k$, $\vr{u_k, z_k}\ge 0$. Therefore, each summand of the right hand side of (\ref{x2}) is non-negative. So we have 
\begin{align}\label{x3}
\sum_{i=1}^m \lambda_i-r\ge \lambda_k (1+\vr{u_k, z_k}).
\end{align}
Summing up for $k$, we have $m(\sum_{i=1}^m \lambda_i-r) \ge \sum_k\lambda_{k=1}^m (1+\vr{ u_k, z_k})\ge a(\sum_{i=1}^m \lambda_i-r)$. By (\ref{x3}), $\sum_{i=1}^m \lambda_i-r>0$, so we have $a \le m$, as desired.

{\bfseries Conditions on the equality.} Suppose $\rho_X=2d$. Then all the inequalities (\ref{3parts}-\ref{x3}) above are equalities, and $\vv{B_{+}}=a=m$.
For any fixed $k$, $1\le k\le m$, The equality of (\ref{x1}) means that $(B_k, G_k$) is also a minimizing pair. By (\ref{x3}), $\delta_k:=\vr{u_k, z_k}=(\sum_{i=1}^m \lambda_i-r)/\lambda_k-1$ is independent of the choice of $z_k$ and $B$. This with (\ref{x1}) implies that
\[
 \vr{u_k, w}= m-1-\delta_k
\]
is also independent of $z_k$ and $B$. By Proposition~\ref{zeroCovering} (3), any neighboring facet $G\in \adjt(F)$ can be expanded into $(G=G'_1,\cdots G'_{m'})$ satisfying Proposition \ref{zeroCovering} (1). Running the same argument above for $(G'_1,\cdots G'_{m'})$, we find $m'=m$, and $\vr{u_G, w}$ is a constant depending on $G$. Since $u_F$ and $\{u_G: G\in \adjt(F)\}$ together span $M_\R$, and $\vr{u_F, w}=m$ is constant, these values decide the unique vector $w\in N$, which only depends on $F$. Consider $G=G_1'$, and the minimizing numerical basis $B_1'$ constructed as above by adding $z_1'$. Then for the minimizing pair $(B_1', G_1')$, $w(B_1', G_1')$ is unique regardless of the choice of $z_1'$. Since
\begin{align}\label{x4}
w(B_1', G_1')=z_1' + \sum_{b\in (B_{0}\sqcup B_{+})\backslash G_1'} b,
\end{align}
we conclude that there is a unique $z_1'$. Since this holds for any $G\in \adjt(F)$, $F$ is a simplex.
We have shown that all the neighboring facets of $F$ are minimizing a numerical basis too. Inductively, we conclude that all the facets of $P$ are simplices. In other words, $P$ is simplicial. By \cite{Casagrande2005}, $d$ is even and $P$ is a free sum of $d/2$-many copies of the reflexive hexagon $\mathcal{H}$.\qed

The numerical bases, and the technique of varying one basis to obtain another, are worth further investigation. A structure theory of the pairs $(B,F)$ could reveal more information on the toric Picard groups. It would also be very interesting to find a more algebro-geometric proof of our main result.

\subsection{Analyzing the proof found using ChatGPT \gptversion{}} 

All steps of this AI-proof are generalizing existing ideas and proofs. For the covering of the zero layer, Proposition \ref{zeroCovering} improves \cite[Lemma 5.5]{Nill2005M}, which is used in Casagrande's proof \cite{Casagrande2005} and also in \cite{Fujita2019}. The AI used the refined Steinitz’s Theorem (Lemma \ref{SteinitzRefined}) without referring explicitly to \cite[Lemma 3.4]{CJT2026}, while the original Steinitz’s theorem has already been used in \cite{Debarre2000} to prove the $2d^2-d$ upper bound by \cite{Voskresenskiĭ1985}, and later in \cite[Thm. 5.4]{Nill2005M}. 

\label{AI-subsec}
The bound on the positive part also combines multiple established techniques. The strategy of taking the facet $F$ minimizing $F(B)$ is introduced in \cite{2008Øbro}, and the spanning by $S$ of the quotient space dual to $Z_F$ is used in \cite{Casagrande2005}. Both proofs arrive at the estimation of $\vv{B_{+}}\le \vv{B_{-}}$ along the direction $u_F$. Eikelberg's theory, and the later reintroduction by \cite{Fujita2019} and \cite{He2026} are essential for bounding simultaneously the zero layer and the positive part. In particular, the concept of numerical basis, and Lemma \ref{negativeCount} are first introduced by \cite{Fujita2019}, while Corollary \ref{1perstep} first appeared in \cite[Lemma 2.16]{Eikelberg1993}, and is the core of the proof of \cite[Thm 3.19]{He2026} controlling the dimension of $H_\Q/\AD(X)$ via chains of facets where each raises the accumulated dimension by at most one.

\bibliography{mybib.bib}
\bibliographystyle{alpha}
\end{document}